\documentclass[12pt]{amsart}

\usepackage{amssymb}
\newtheorem{theorem}{Theorem}
\newtheorem{lemma}{Lemma}

\newcommand{\on}{\operatorname}
\newcommand{\Eul}{\on{Eul}}
\newcommand{\Res}{\on{Res}}
\newcommand{\Vol}{\on{Vol}}
\newcommand{\SU}{\on{SU}}
\newcommand{\Ad}{\on{Ad}}
\newcommand{\su}{\on{su}}
\newcommand{\U}{\on{U}}

\newcommand{\ra}{\rightarrow}
\newcommand{\F}{\mathcal{F}}
\newcommand{\f}{\frac}
\newcommand\om{\omega}
\newcommand{\Z}{\mathbb{Z}}

\newcommand{\C}{\mathbb{C}}
\newcommand\lie[1]{\mathfrak{#1}}
\newcommand{\g}{\lie{g}}
\renewcommand{\t}{\lie{t}}
\renewcommand{\d}{\on{d}}
\newcommand{\ti}{\tilde}
\newcommand{\ot}{\overline{\theta}}

\begin{document}

\author{Olga Plamenevskaya}
\address{Department of Mathematics \\ St.Petersburg University \\
  \ St.Petersburg \\ Russia}
\email{Olga.Plamenevskaja@pobox.spbu.ru}

\title{A residue formula for $\SU(2)$-valued moment maps }

\date{\today}

\begin{abstract}
In ~\cite{JK} Jeffrey and Kirwan suggested expressions
for intersection pairings on the reduced space
$M_0=\mu^{-1}(0)/G$ of a Hamiltonian $G$-space $M$
in terms of multiple residues.
In this paper we prove a residue formula for
symplectic volumes of reduced spaces of a quasi-Hamiltonian
$\SU(2)$-space. The definition of quasi-Hamiltonian
$G$-spaces was recently introduced in ~\cite{AMM}.
\end{abstract}

\maketitle

\section{Introduction}

Let $(M, \omega)$ be a compact symplectic manifold with
Hamiltonian action of a compact Lie group $G$ and $\mu$
be the corresponding moment map.
Jeffrey and Kirwan  ~\cite{JK} gave expressions
for general intersection pairings on the homotopy quotient
$M_0=\mu^{-1}(0)/G$ in terms of multiple residues of
certain integrals over connected components of the fixed point
set $ M^T $ of the maximal torus $ T \subset G$.

In~\cite{AMM} and~\cite{AMW}, Alekseev, Malkin, Meinrenken  and Woodward
 develop a theory of quasi-Hamiltonian
$G$-spaces with the moment map taking values in a Lie group $G$. For
such spaces, they
introduce analogs of "classical" Hamiltonian reduction, Liouville
volumes, Duistermaat-Heckmann measure, localization formulas, etc.
The goal of this note is to obtain a  quasi-Hamiltonian
residue localization formula for the special case $G=\SU(2)$.

We first introduce
some notation. Let  $T=S^1\subset\SU(2)$ be a choice
of a Cartan circle of $G=\SU(2)$, $\t,\g$ the corresponding Lie algebras.
Fix a positive Weyl chamber $\t_+\subset\t$. Let $\alpha$ be the
unique positive root of $\t$, and $\rho\in\t^*$ the fundamental
weight. We identify $\g$ and $\g^*$ via a scalar product on $\g$
such that $(\alpha,\alpha)=2$, and write $\alpha=2\rho$.
Note that this choice of the inner product implies
that $\Vol T=\sqrt{2}$, $\Vol G=\f{\sqrt{2}}{2\pi}$.
We parameterize the alcove $\mathfrak{A}\cong
[0,1]$,  $t\in [0,1]$ corresponding to $t\rho\in\t$.

We shall briefly describe the results of~\cite{AMM} and~\cite{AMW}.
As we are concerned with the case $G=\SU(2)$, we give below
 the formulas applied
to $\SU(2)$ and do not state the results for the general case of a compact
connected group $G$.
Let $\theta=g^{-1}\d\!g$ and $\ot=\d \!g\/ g^{-1}$ denote the left- and right-invariant
Maurer-Cartan forms on $G=\SU(2)$, and let $\chi$ be the canonical closed bi-invariant
$3$-form on $G$,
$$
\chi=\f{1}{12}(\theta,[\theta,\theta])=\f{1}{12}(\ot,[\ot,\ot]).
$$

{\bf Definition.} \cite{AMM} A quasi-Hamiltonian $G$-space is a
$G$-manifold
with a invariant $2$-form $\om\in\Omega(M)^{G}$ and an equivariant map
$\Phi\in C^{\infty}(M,G)^G$ (the "moment map"), such that:
\begin{enumerate}
\item[(i)]
The differential of $\omega$ is given by
$\d\om=-\Phi^{*}\chi$.
\item[(ii)]
The moment map satisfies
$$
\imath(v_\xi)=\f{1}{2}\Phi^{*}(\theta+\ot,\xi).
$$
\item[(iii)]
At each $x\in M$, the kernel of $\om_x$ is given by
$$
\ker \om_x =\{ v_\xi, \xi\in\ker(\Ad_{\Phi(x)}+1)\}.
$$
\end{enumerate}

Similarly to Meyer-Marsden-Weinstein reduction, the q-Hamiltonian
reduced phase spaces are defined in the following way.
Let $g\in M$ be a regular value of the moment map $\Phi$, then the pre-image
$\Phi^{-1}(g)$ is a smooth submanifold on which the action of the centralizer
$G_g$ is locally free. The reduced
space $M_g=\Phi^{-1}(g)/G_g$ is a symplectic orbifold.

As in the Hamiltonian case, there is the volume form $\Gamma$ on $M$,
called the Liouville  form.
For $G=\SU(2)$ it is defined as follows.
By equivariance of the moment map it suffices to define $\Gamma$
for points $x\in M$ such that $\Phi(x)\in T$.
 If $\Phi(x)=\exp(t\rho)$, and
$t\neq \pm 1/2$, then
$$
\Gamma_{x}=\f{1}{\cos \pi t}(\exp \omega)_{[top]}.
$$
The denominator $\cos \pi t$ cancels the zeroes of $\exp \om$
arising because of condition~(iii), so $\Gamma$  extends smoothly
to all points of $\Phi^{-1}(T)$.

The Duistermaat-Heckman measure $\varrho$ on $G$ is
now defined as the push-forward of the Liouville form under the moment map,
$$
\varrho=\Phi_*\Gamma.
$$
We shall consider the DH function, i.e., the density of
this measure with respect to the Haar
measure, and use the same notation for it, writing
$$
\varrho=\varrho(g)\d \Vol_G.
$$
As in the Hamiltonian setting, the Duistermaat-Heckmann measure is related
to volumes of the reduced spaces. Let $g\in G$ be a regular
value of $\Phi$, and $M_g$ the corresponding reduced space.
Then $\varrho$ is smooth at $g$.
Let $k$  be the cardinality of a generic stabilizer for $G_g$-action
on $\Phi^{-1}(g)$.
Then for the reduced space at
$g=\exp(t\rho)$ with $t\in (0, 1)$ we have
\begin{equation}\label{vol}
\Vol(M_g)=k \f {2 \sin \pi t}{\sqrt{2}}{\varrho(g)},
\end{equation}
and for central elements $g=\pm e$
\begin{equation}\label{volc}
\Vol(M_g)=k\f{2\pi}{\sqrt{2}}{\varrho(g)}.
\end{equation}

In the case of Hamiltonian torus actions, the localization formula
by Berline-Vergne~\cite{BV} gives for the Fourier-Laplace transform
of the DH function in terms
of integrals over
fixed point manifolds of the torus  subgroup generated by $\xi$,
$$
\int_M e^{i \pi \langle \Phi, \xi \rangle} \exp\om=
\sum_{F\in \F(\xi)}\int_{F}\f{\exp(\om_F)}
{\Eul(\nu_F,\xi)}e^{i \pi \langle\Phi,\xi\rangle},
$$
where $\Eul(\nu_F,\xi)$ is the equivariant Euler class.
The q-Hamiltonian counterpart of this formula, valid for arbitrary compact connected
Lie group $G$, deals with
the Fourier coefficients of the Duistermaat-Heckmann
function, which localize on the fixed point
manifolds of certain circle subgroups of $G$.

Let $\F$ be the set of connected components $F\subset M$ of the
fixed point set of the Cartan circle $T\subset\SU(2)$.
Each $F\in\F$ is a symplectic manifold with the pull-back
$\om_F$ of $\om$ as a $2$-form.
By equivariance, the restriction $\Phi|_{F}$
is  constant and sends $F$ to a point  of $T$.
Writing
$\Phi_F=\exp(\mu_F\rho)$,
we introduce a real function $\mu$ with values in $(-1,1]$.
Let $\Phi_F^{n\rho}$ stand for $e^{\pi i n\mu_F}$.
We denote by$\F_+$ the set of components $F\in\F$ with
$\Phi\in \mathfrak{A}$, that is, $\mu_F\geq 0$.
Orientations of $M$ and $F$ induce the orientation on the normal
bundle $\nu_F$, and the $T$-equivariant Euler
class $\Eul(\nu_F,\cdot)$ is defined.

Recall that irreducible representations of $\SU(2)$ are
labelled by their highest weights $n\rho$, $n=0,1,2,\dots$.
The  dimension of the representation  $V_n$ with
the highest weight $n\rho$ is $\dim V_n =n+1$.
Let $\chi_n$ denote the character of $V_n$.

\begin{theorem}[q-Hamiltonian localization formula,~\cite{AMW}]\label{loc}
The Fourier coefficients of the DH function are given by
\begin{equation}\label{CoeffsLoc}
\langle\varrho,\chi_n\rangle= \dim V_n\, \sum_{F\in\F}
\int_F
\f{\exp(\om_F)}{\Eul(\nu_F,2\pi i(n+1)\rho)}
\,\Phi_F^{(n+1)\rho}.
\end{equation}
\end{theorem}
The DH function is then reconstructed as
$\varrho(g)=\f{1}{\Vol G}\sum_n \langle
\varrho,\chi_n \rangle \chi_n(g^{-1})$.

\section{The residue formula}

The purpose of this note is to obtain  a quasi-Hamiltonian
residue formula for the Liouville
volumes of the reduced phase spaces for
$\SU(2)$-actions. Because of~(\ref{vol}) and~(\ref{volc}),
this is the same as to get a formula for the DH function $\varrho(g)$.
Since  $\varrho(g)$ is a function of conjugacy classes, it suffices
to evaluate it
at the elements of the Cartan circle, $g=\exp(t\rho)$.

\begin{theorem}\label{residue} The function $\varrho(t)$ is a sum
of contributions $\varrho_F(t)$  of the components $F\in\F_+$,
$$
\varrho(\exp(t\rho))=\sum_{F\in\F_+}\varrho_F(\exp(t\rho)),
$$
where $\varrho_F(\exp(t\rho))$  corresponding to $\mu_F\in(0,1)$ are given by
\begin{equation}\label{t<mu}
\varrho_F(\exp(t\rho))=-\f{4 \pi^2 i}{\sqrt{2}}\,{\sin\pi t}\Res_0
\f{ze^{\pi i z\mu_F}\sin (\pi t z)}{e^{2\pi iz}-1}
\int_F \f{\exp(\om_F)}{\Eul(\nu_F,2\pi i z\rho)}
\end{equation}
for $0<t<\mu_F$, and
\begin{equation}\label{t>mu}
\varrho_F(\exp(t\rho))=\f{4 \pi^2 i}{\sqrt{2}}\,{\sin \pi t}\Res_0
\f{ze^{\pi i z(\mu_F+1)}\sin (z\pi(1-t))}{e^{2\pi iz}-1}
\int_F \f{\exp(\om_F)}{\Eul(\nu_F,2 \pi i z\rho)}
\end{equation}
for  $\mu_F<t<1$.

If $\mu_F=0$ or $\mu_F=1$, the formulas~(\ref{t<mu}) and~(\ref{t>mu})
are valid if we add the factor $1/2$.

If  $t=0$ or $t=1$ is a regular value of $\varrho(\exp(t\rho))$ (i.e.,
$e$ or $-e$ is a regular value of the moment map), then
\begin{eqnarray}
\varrho_F(e)&=-&\f{4 \pi^2 i}{\sqrt{2}}\,\Res_0
\f{z^2 e^{\pi i z\mu_F}}{e^{2\pi iz}-1}
\int_F \f{\exp(\om_F)}{\Eul(\nu_F,2\pi i z\rho)}\label{0}, \\
\varrho_F(-e)&=&\f{4 \pi^2 i}{\sqrt{2}}\,\Res_0
\f{z^2 e^{iz(\mu_F+\pi)}}{e^{2\pi iz}-1}
\int_F \f{\exp(\om_F)}{\Eul(\nu_F,2 \pi i z\rho)}.\label{1}
\end{eqnarray}
\end{theorem}

Note that in the case $G=\SU(2)$ the Jeffrey--Kirwan formula for
the Liouville volume of $\mu^{-1}(0)/G$ (i.e., the value of the standard
DH function at zero) gives
$$
\varrho(0)=-\f{1}{2}\Res_0 z^2\sum_{F\in\F_+}e^{iz\mu_T(F)}\int_F
\f{e^{i\om_F}}{\Eul(\nu_F,z)},
$$
$\F_+$ denoting the set of the components $F$ for which $\mu_T(F)>0$.
Even in the case of Hamiltonian $G$-spaces
our formulas look different as we work with discrete Fourier series instead
of continuous Fourier transform used in~\cite{JK}.

In the proof of the Theorem, we shall use the following lemma.
\begin{lemma}
 Let $f(z)$ be a rational function with the only pole at zero,
 such that $f(z)\rightarrow 0$ as $z\rightarrow\infty$. Then for
 $0<\gamma<2\pi$

\begin{equation}\label{Szenes}
  \sum_{m\in\Z,m\ne 0} e^{i m\gamma}f(m)=-2\pi i\Res_0\f{f(z)e^{i\gamma z}}
  {e^{2\pi i  z}-1},
\end{equation}

 and for $-2\pi<\gamma<0$
\begin{equation}\label{Szenes'}
\sum_{m\in\Z,m\neq 0} e^{im\gamma} f(m)=
-2\pi i\Res_0 \f{f(z)e^{i\gamma z}}{1-e^{-2\pi i z}}.
\end{equation}

\end{lemma}
\begin{proof} The proof of the Lemma is a straightforward calculation,
which is given in the Appendix.
\end{proof}

Note that this Lemma resembles to a one-dimensional case of results by
Szenes~\cite{S} and Jeffrey-Kirwan\cite{JK*}.

\begin{proof}[Proof of the Theorem~\ref{residue}]
We shall make use of Theorem~\ref{loc} and sum up
the Fourier series. Introduce the functions
$\tilde{\varrho}_F(\exp(t\rho))$ with the Fourier coefficients
\begin{equation}\label{CoeffF}
\langle\tilde{\varrho}_F,\chi_n\rangle=\dim V_n \,
\int_F \f{\exp(\om_F)}{\Eul(\nu_F,2 \pi i (n+1)\rho)}
\, \Phi_F^{(n+1)\rho},
\end{equation}
so that by~(\ref{CoeffsLoc})
$$
\varrho(t)=\sum_{F\in\F}\ti{\varrho}_F(t).
$$

Recall that all the characters of $\SU(2)$ are real and have the
property $\chi_n(g)=\chi_n(g^{-1})$, so we can just write
$\varrho(g)=\f{2\pi}{\sqrt{2}}\sum\langle\varrho,\chi_n\rangle\chi_n(g)$. Further, the
character $\chi_n$ is given at $g=\exp(t\rho)$ by the Weyl
character formula,
$$
\chi_n(\exp(t\rho))=\frac{\sin(\pi(n+1)t)}{\sin\pi t}.
$$
Substituting this into~(\ref{CoeffF}) and
putting together the Fourier series, we get
$$
\tilde{\varrho}_F(\exp(t\rho))=\f{2\pi}{\sqrt{2}}\sum_{n\geq 0}\f{n+1}{2 i \sin\pi t}
(e^{\pi i (n+1)t}-e^{\pi i (n+1)t})
\int_F \f{\exp(\om_F)}{\Eul(\nu_F, 2\pi i(n+1)\rho)}
\,\Phi_F^{(n+1)\rho}.
$$

We now gather the components $F\in\F$ of
the fixed point set into pairs $F, F'=F^w$ with
$F\in\F_+$ by means of the Weyl element
$w=\mbox{\tiny$\left( \begin{array}{cc} 0 & 1 \\ -1 & 0 \end{array}\right)$}$.
For $F\in\F$  with $\mu_F\neq 0, 1$ we let
$\varrho_F=\ti{\varrho}_F+\ti{\varrho}_{F'}$. If $\mu_F=0$ or $1$,
then $F=F'$, and we denote $\varrho_F=\ti{\varrho}_F
=\f{1}{2}(\ti{\varrho}_F+\ti{\varrho}_{F'})$. Then
$$
\varrho(\exp(t\rho))=\sum_{F\in\F_+}\varrho_F(\exp(t\rho)).
$$
 Use the change of variables $y=x^w$
to relate the functions $\ti{\varrho}_F$ and $\ti{\varrho}_{F'}$.
For $h\in S^1$ we have $hw=wh^{-1}$,
so the action of $2\pi i n\rho$ at $\nu_{F'}$ corresponds to
the action of $-2 \pi i n\rho$ at $\nu_F$.
The form $\om$ is invariant, then
$$
\int_{F'} \f{\exp(\om_{F'})}{\Eul(\nu_{F'},2 \pi i n\rho)}=
\int_F \f{\exp(\om_F)}{\Eul(\nu_F,-2\pi i n\rho)}.
$$
By the equivariance of the moment map,
$\Phi_{F'}=\Phi_F^{-1}$.

Now pick together the terms with $e^{\pi i(n+1)t}$ from
$\ti{\varrho}_F$ and those with $e^{-\pi i(n+1)t}$
from $\ti{\varrho}_{F'}$, and vice versa. Substituting $m=n+1$ and
$\Phi_F^{n\rho}=e^{\pi i n\mu_F}$, we get
\begin{equation}\label{mSeries}
\varrho_F(\exp(t\rho))=\f{2\pi}{\sqrt{2}}\f 1{2i\sin{\pi t}}\sum_{m\in\Z, m\neq 0}
m (e^{\pi i m (t+\mu_F)}-e^{-\pi i m (t-\mu_F)})
\int_F \f{\exp(\om_F)}{\Eul(\nu_F,2 \pi i m\rho)}.
\end{equation}
(Here we assume that $\Phi_F\neq\pm e$; otherwise a factor $1/2$
must be added).

We assume that the moment map has regular points. It follows
that the codimension
of $F$ in $M$ must  be at least $4$. Then the rational function
$1/{\Eul(\nu_F, 2\pi i z\rho)}$ is of
order of at least $1/z^2$.
It means that (\ref{mSeries}) consists of
two series to which the Lemma applies.

If $t<\mu_F$, both series sum up by means of~(\ref{Szenes}),
 and we get
$$
\varrho_F(t)=-\f{4\pi^2 i}{\sqrt{2}}{\sin \pi t}\Res_0
\f{z e^{\pi i z\mu_F}\sin (\pi t z)}{e^{2\pi i z}-1}
\int_F \f{\exp(\om_F)}{\Eul(\nu_F,2 \pi i z\rho)},
$$
which is exactly~(\ref{t<mu}).
If $t>\mu_F$, we  use both~(\ref{Szenes}) and~(\ref{Szenes'}),
and obtain the formula~(\ref{t>mu}).
To get the formulas for $\varrho_{F}$ with $\Phi_F=\pm e$ we must
just add a factor $1/2$ into~(\ref{t<mu})and~(\ref{t>mu}).

If $e$ or $-e$ are the regular values of the moment map, we
can consider the limits $t\ra 0$ and $t\ra 1$ to
obtain the formulas~(\ref{0}) and~(\ref{1}), respectively.
\end{proof}

{\bf Example 1.} We apply the obtained residue formulas to the following
example of a q-Hamiltonian $\SU(2)$-space from~\cite{AMW}.
The space is constructed as follows. Take $\C^2$ equipped with
its natural symplectic form and defining $\SU(2)$-action, and
let $\Phi_0:\C^2\rightarrow\su(2)$ be its classical moment map.
Let $Y_1=Y_2\subset \C^2$ be the open ball given as the pre-image
$\Phi_0^{-1}(G\cdot[0,\rho))$, and $\Phi_{1,0}=\Phi_{2,0}$ the
restrictions of $\Phi_0$. Then
$$
Y_3=\Phi_{1,0}^{-1}((0,\rho))
$$
as a Hamiltonian $\U(1)$-space is equivariantly symplectomorphic
to
$$
Y'_3=\Phi_{2,0}^{-1}((-\rho,0))
$$
via the isomorphism $(z,\xi)\rightarrow(z,\xi-\rho)$.
We now glue the spaces $Y_1$ and $Y_2$ together along their
boundaries by means of
the embeddings
$$
Y_1\leftarrow \SU(2)\times_{U(1)}Y_{3}\rightarrow Y_2,
$$
and obtain a sphere $S^4$ with $\SU(2)$ acting by rotations.

The action has two fixed points, one with
$\Phi=e$, the other
with $\Phi=-e$.The Euler classes
are given by $\Eul(\nu_F,\xi)=\mp\langle\rho,\xi\rangle^2$
for $\xi\in \t$, so $\Eul(\nu_F, 2\pi i z\rho)=\pm\pi^2 z^2$.

Denoting by $\varrho_0$, $\varrho_1$
the contributions of the fixed points, we write
$$
\varrho((\exp(t\rho))=\varrho_0((\exp(t\rho))+\varrho_1((\exp(t\rho)).
$$
Using the formula~(\ref{t<mu}) (with factor 1/2), we get
$$
\varrho_1((\exp(t\rho))=\frac{2 \pi^2 i}{\sqrt{2}}{\sin \pi t} \Res_0
\frac{z e^{\pi i z t}\sin(\pi t z)}{e^{2\pi i z}-1}
\frac{1}{\Eul(\nu_1,2 \pi z\rho)}.
$$
Computing the residue, we get
$$
\varrho_1((\exp(t\rho))=\frac{t}{\sqrt{2} \sin \pi t}.
$$
In the same way
$$
\varrho_0((\exp(t\rho))=\frac{1-t}{\sqrt{2} \sin \pi t},
$$
hence
$$
\varrho((\exp(t\rho))=\frac{1}{\sqrt{2}\sin \pi t}.
$$
Now use~(\ref{vol}), and get $\Vol M_g=1$. This is obviously the
correct answer, because the reduced spaces are just points.

{\bf Example 2.}
Let us now apply our residue formula to an
$\SU(2)$-space $\SU(2)^{2n}$ obtained by means of a fusion product
introduced in~\cite{AMM}.

Given a q-Hamiltonian $G\times G$-space $M$ with  $2$-form moment map
$\Phi=(\Phi_1,\Phi_2)$, we can consider $M$ as a $G$-space
with a diagonal $G$-action. Then $M$ with moment map
$\ti \Phi=\Phi_1\Phi_2$ and $2$-form $\ti \om=\om+\f{1}{2}(\Phi_1^*\theta,
\Phi_2^*\ot)$ is a q-Hamiltonian $G$-space.
Taking two q-Hamiltonian $G$-spaces $M_1$, $M_2$, we
apply the above procedure to $G\times G$-space $M_1\times M_2$,
and obtain the fusion product $M_1\circledast M_2$.

Now, following~\cite{AMM}, consider a
double  $D(\SU(2))$, that is, the $SU(2)\times \SU(2)$-space
$\SU(2)\times\SU(2)$ defined as follows. The group action is given by
$$
(a,b)^{(g_1,g_2)}=(g_1 a g_2^{-1}, g_2 b g_1^{-1}),
$$
the moment map is $\Phi=(\Phi_1,\Phi_2)$, where
$$
\Phi_1(a,b)=ab,\quad \Phi_2(a,b)=a^{-1}b^{-1},
$$
and the $2$-form is
$$
\om_D=\f{1}{2}(a^*\theta, b^*\ot)+\f{1}{2}(a^*\ot, b^*\theta).
$$
We can get a $G$-space ${\bold D}(\SU(2))$, applying fusion to
$D(\SU(2))$. The group $\SU(2)$ acts by conjugations on each factor
of ${\bold D}(\SU(2))=\SU(2)\times\SU(2)$, and the moment map is
$\Phi(a,b)=aba^{-1}b^{-1}=[a,b]$.

We shall consider the q-Hamiltonian $\SU(2)$-space
$\SU(2)^{2n}$ obtained as a fusion  product of $n$
copies of ${\bold D} (\SU(2))$.
The fixed point set of the $T$-action has only one component
$F$, which is a product of $2n$
copies of $T\subset\SU(2)$. The moment map $\Phi$ sends this set to $e$,
so we must insert the factor $1/2$ into the formula~(\ref{t>mu}).
To compute the Euler class $\Eul(\nu_F,\cdot)$, note that the normal bundle
is trivial and can be identified with $(\g/\t)^{2n}$, and
$\Eul(\nu_F,2 \pi z\rho)=z^{2n}(\langle 2\rho,2\pi i \rho\rangle
\langle -2\rho, 2\pi i \rho\rangle)^{n}=(2z)^{2n}\pi^{2n}$.

For the double $D(\SU(2))$, the form $\exp(\ti \om_{T\times T})$
coincides with the Riemannian volume form induced
by our choice of the inner product on $\g$. As the restriction of the moment
map to the torus in the fusion product is trivial, at points of
the torus this form remains intact under fusion.
Then the form $\omega_F$ on  $F=T^{2n}$ is just the Riemannian volume form,
and we have
$$
\int\exp(\omega_F)=(\Vol T)^{2n}=2^n.
$$

Substituting everything into~(\ref{t>mu}),
we get
$$
\varrho(\exp(t\rho))=\f{2\pi}{\sqrt{2}}\f{i}{2^n \pi^{2n-1} \sin \pi t}
\Res_0 \f{e^{\pi i z}\sin(z\pi(1-t))}{z^{2n-1}(e^{2\pi i z}-1)}.
$$
This can be re-written as
$$
\varrho(\exp(t\rho))=\f{i\sqrt{2}}{2^{n}\pi ^{2n-2}(2n-2)!\sin{\pi t}}
\f{\partial^{2n-2}}{\partial z^{2n-2}}\left.\left(\f{e^{\pi i z}\sin(z\pi(1-t))}{e^{2\pi i z}-1}\right)\right|_{z=0}.
$$
If $n=1$,  for  the space ${\mathbf D}(\SU(2))$ we get
$$
\varrho(\exp(t\rho))=\f{1-t}{2\sqrt{2}\sin\pi t}
$$
for $t\in(0,1)$.

To find the volumes of the reduced spaces, note that the generic
stabilizer is just the center of $\SU(2)$, which consists of two points.
Using~(\ref{vol}), we have
$$
\Vol(M_g)=1-t.
$$
This coincides with the answer given by Witten's
formula~{\cite{W}} and~{\cite{JK*}}  for the volumes of moduli spaces.

\appendix
\section{}
We now give the proof of the Lemma.
\begin{proof} We have $f=a_1/z+a_2/{z^2}+\dots+a_n/{z^n}$, so it
suffices to prove~(\ref{Szenes}) for $f=1/{z^k}$, $k \leq 1$.
Consider a function
$$
h(z)=\f{f(z)e^{i\gamma z}}{e^{2\pi i z-1}}.
$$
Its poles are the real integer numbers, and $\Res_{z=m}h(z)=f(m)/2\pi i $
for $m\in\Z$, $m\ne 0$.  Integrate $h(z)$ over the contour
\begin{align*}
    \Gamma_m=&\{-m-1/2\leq x\leq m+1/2,\ y=m\}\cup \\
             &\{x=m+1/2,\-m\leq y\leq m\}\cup     \\
             &\{-m-1/2\leq x\leq m+1/2,\ y=-m\}\cup \\
             &\{x=-m-1/2,\ -m\leq y\leq m\}.
\end{align*}

By the Cauchy theorem
 $$2\pi i \sum_{-l \leq m\leq l}\Res_{z=m} h(z)=
\int_{\Gamma_l}h(z)\,d z,
$$
and we have only to prove that $\int_{\Gamma_l}h(z)\,d z\rightarrow 0$ as
$m\rightarrow 0$.

While $z=(x,y)\in \Gamma_{1, m}=\{-m-1/2\leq x\leq m+1/2,\ y=m\}$,
$|e^{i\gamma z}/(e^{2\pi i z}-1)|=|e^{i\gamma (x+im)}/(e^{2\pi i(x+im)}-1)|
\leq e^{\gamma m}/(e^{2\pi m}-1)$. Then
$$
\left|\int_{\Gamma_{1, m}}h(z)\,d z\right|\leq
(2m+1)e^{\gamma m}/m (e^{2\pi m}-1),
$$ which tends to $0$ since $\gamma<2\pi$.

For $\Gamma_{2, m}=\{-m-1/2\leq x\leq m+1/2,\ y=-m\}$ we get
$$
\left|\int_{\Gamma_{2, m}}h(z)\,d z\right|\leq
(2m+1)e^{-\gamma m}/m (e^{-2\pi m}-1),
$$ which tends to $0$ because
$\gamma>0$.

Now cut $\Gamma_{3, m}=\{x=m+1/2,\ -m\leq y\leq m\}$ into
\begin{align*}
  \Gamma_{3, m}&=\Gamma'_{3, m}\cup\Gamma''_{3, m}= \\
&=\{x=m+1/2,\ -\sqrt{m}\leq y\leq \sqrt{m}\}
\cup \{x=m+1/2,\  y\geq \sqrt{m}\}.
\end{align*}
If $z=(x,y)\in\Gamma_{3 m}$, then $h(z)=e^{-\gamma y}/(z(e^{-2\pi
y}+1))$. It follows that
$$\left|\int_{\Gamma'_{3 m}}h(z)\,d z\right|
\leq 2\sqrt{m}/m.
$$
If $z\in \Gamma''_{3 m}$, we refine the estimate.
We have $e^{-\gamma y}/(e^{-2\pi y}+1)\leq e^{-\gamma\sqrt{m}}$ for
$y\geq\sqrt{m}$, and
$e^{-\gamma y}/(e^{-2\pi y}+1)\leq e^{-(2\pi-\gamma)\sqrt{m}}$
for $y\leq\-\sqrt{m}$. Therefore
$$
\left|\int_{\Gamma''_{3 m}}h(z)\,d z\right|\leq
e^{-\gamma\sqrt{m}}+e^{-(2\pi-\gamma)\sqrt{m}}.
$$
It follows that
$\int_{\Gamma_{3,m}}h(z)\,d z\rightarrow 0$.

The integral over $\Gamma_{4 m}=\{x=-m-1/2,\ -m\leq y\leq m\}$
is treated in the same way.
\end{proof}

I thank A. Alekseev for his helpful suggestions and encouragement.

\end{document}